\documentstyle[11pt,leqno]{article}
\newtheorem{theorem}{{\sc Theorem}}
\newcommand{\bt}{\begin{theorem}}
\newcommand{\et}{\end{theorem}}
\newcommand{\newsection}[1]{\setcounter{equation}{0} \setcounter{theorem}{0}
\section{#1}}

\newcommand{\NI}{\noindent}
\newcommand{\bea}{\begin{eqnarray}}
\newcommand{\eea}{\end{eqnarray}}

\def \spec#1 {\mathop{#1}}

\def \b #1 {\bf #1}

\newcommand {\CC}{\centerline}

\newcommand{\cle}{{\cal E}}
\newcommand{\clf}{{\cal F}}

\newcommand{\clh}{{\cal H}}

\newcommand{\ity}{\infty}
\newcommand{\raro}{\rightarrow}

\newcommand{\vsp}{\vskip 1em}

\newcommand{\be}{\begin{equation}}
\newcommand{\ee}{\end{equation}}
\newcommand{\ben}{\begin{eqnarray*}}
\newcommand{\een}{\end{eqnarray*}}

\begin{document}
\CC{\bf{Nonparametric Estimation of Linear Multiplier for} }
\CC{\bf{Processes Driven  by a bifractional Brownian Motion}}
\vsp
\CC{B.L.S. Prakasa Rao}
\CC{CR RAO Advanced  Institute of Mathematics, Statistics}
\CC{and Computer Science, Hyderabad, India}
\CC{(e-mail address: blsprao@gmail.com)}
\vsp
\NI{\bf Abstract:} We study the problem of nonparametric estimation of the linear multiplier function $\theta(t)$ for processes satisfying stochastic differential equations of the type
$$dX_t= \theta(t)X_t dt+ \epsilon\;dW_t^{H,K}, X_0=x_0, 0 \leq t \leq T$$
where $\{W_t^{H,K}, t\geq 0\}$ is a bifractional Brownian motion with known parameters $H\in (0,1), K\in (0,1]$ and $HK\in (\frac{1}{2},1).$We study the asymptotic behaviour of the estimator of the unknown function $\theta(t)$ as $\epsilon \raro 0.$
\vsp
\NI{\bf Keywords :} Nonparametric estimation, Linear multiplier, bifractional Brownian motion.
\vsp
\NI{\bf Mathematics Subject Classification :} Primary 60G22, Secondary 62G05.
\newsection {Introduction}

Statistical inference  for fractional diffusion processes satisfying stochastic
differential equations driven by a fractional Brownian motion (fBm)  has been studied
earlier and a comprehensive survey of various methods is given in Mishura (2008) and Prakasa
Rao (2010). fBm is the only self-similar Gaussian process with stationary increments starting from zero. For small increments, in models for turbulence, fBm is considered as a good model but inadequate for modeling large increments. Houdr\'e and Villa (2003) introduced a process called a bifractional Brownian motion that can be considered as a generalization of a fBm which retained the properties of self-similarity, stationarity for small increments and increased the choice of processes for modeling phenomenon such as turbulence. This process is also a quasi-helix as defined in Kahane (1981, 1985). There has been a recent interest to study problems of statistical inference for
stochastic processes driven by a bifractional Brownian motion (bifBm). Keddi et al. (2020) investigated the problem of nonparametric estimation of the trend for processes driven by a bifractional Brownian motion following the methods in Kutoyants (2012), Mishra and Prakasa Rao (2011). Some maximal and integral inequalities for a bifBm were derived in Prakasa Rao (2024).  
\vsp
We now discuss the problem of estimating the function $\theta(t), 0 \leq t \leq T$ (linear multiplier) based on the observations of a  process $\{X_t, 0\leq t \leq T\}$ satisfying the stochastic differential equation
$$
dX_t=\theta(t)\;X_t dt + \epsilon \; dW_t^{H,K}, X_0=x_0, 0 \leq t \leq T
$$
where $\{W_t^{H,K}, t \geq 0\}$ is bifBm and study the properties of the estimator as $\epsilon \rightarrow 0.$ 
\vsp
\newsection {bi-fractional Brownian motion}
We will now describe some properties of a bifractional Brownian motion and properties of processes driven by a bifractional Brownian motion.
\vsp
Let $(\Omega, \clf, (\clf_t), P) $ be a stochastic basis satisfying the
usual conditions and the processes discussed in the following are
$(\clf_t)$-adapted. Further the natural filtration of a process is
understood as the $P$-completion of the filtration generated by this
process.
Consider a centered Gaussian process $W^{H,K}= \{W_t^{H,K}, t\geq 0\}$ called  the {\it bifractional Brownian motion} (bifBm) with the covariance function
$$R_{H,K}(s,t)= \frac{1}{2^K}[(t^{2H}+s^{2H})^K-|s-t|^{2HK}], t\geq 0, s\geq 0$$
where $0<H<1$ and  $0< K\leq 1.$ If $K=1$, then the bifractional Brownian motion reduces to the fractional Brownian motion and if $K=1$ and $H=\frac{1}{2},$ then it reduces to the Brownian motion. As mentioned earlier, a bifBm can be considered as a generalization of the fBm but its increments are not stationary. Russo and Tudor (2006) studied the properties of a bifbm $W^{H,K}.$ Houdr\'e and Villa (2003) and Tudor and Xiao (2007) discussed the following properties of a bifBm $W^{H,K}$ (cf. Tudor (2013, 2023)).
\vsp
\noindent (1) $E(W_t^{H,K})=0, Var(W_t^{H,K})= t^{2HK}, t\geq 0.$ \\
(2) The process $W^{H,K}$ is self-similar with index $HK \in (0,1)$, that is, for every real $a>0,$
$$\{W_{at}^{H,K},t\geq 0\}\stackrel {\Delta}=\{ a^{HK} W_t^{H,K},t \geq 0\}.$$
Here $\Delta$ indicates that the processes, on both sides of the equality sign, have the same finite-dimensional distributions.\\
(3) The process $W^{H,K}$ is not Markov and it is not a semimartingale if $HK \neq \frac{1}{2}.$\\
(4) The sample paths of the process $W^{H,K}$ are Holder continuous of order $\delta$ for any $0 < \delta<HK,$ and they are nowhere differentiable.\\
(5) The bifBm $W^{H,K}$ satisfies the inequalities
$$2^{-K}|t-s|^{2HK}\leq E[W_t^{H,K}-W_s^{H,K}]^2\leq 2^{1-K}|t-s|^{2HK}, t\geq 0, s \geq 0.$$
\vsp
The definition of a bifBm $W^{H,K}$ can be extended for $K\in (1,2)$ with $H \in (0,1)$ and $HK \in (0,1)$ (cf. Bardina and Es-Sebaiy (2011) and Lifshits and Volkava (2015)). Hereafter, we assume that $HK\in (\frac{1}{2},1).$ The  stochastic calculus with respect to the bifractional Brownian motion is developed by Kruk et al. (2007). 
\vsp
Fix a time interval $[0,T]$ and let $\cle $ be the class of real-valued step functions defined on the interval $[0,T].$ Let $\clh_{W^{H,K}}$ be the canonical Hilbert space associated with the bifBm defined as the closure of the set $\cle$ with respect to the inner product
$$<I_{[0,t]},I_{[0,s]}>_{\clh_{W^{H,K}}}=R_{H,K}(t,s)=\int_0^T\int_0^T I_{[0,t]}(u)I_{[0,s]}(v)\frac{\partial^2 R_{H,K}(u,v)}{\partial u\partial v}dudv$$
where $R_{H,K}(t,s)$ is the covariance of the random variables $W_t^{H,K}$ and $W_s^{H,K}.$ The mapping $\varphi:\cle \raro W^{H,K}(\varphi)$ is an isometry from $\cle$ to the Gaussian space generated by $W^{H,K}$ and it can be extended to the space $\clh_{W^{H,K}}$. We consider the subspace $|\clh_{W^{H,K}}|$ of $\clh_{W^{H,K}}$ as the set of measurable functions $\varphi$ on $[0,T]$ satisfying
\be
||\varphi||_{|\clh_{W^{H,K}}|}=\int_0^T\int_0^T \varphi(u) \varphi(v)\frac{\partial^2 R_{H,K}(u,v)}{\partial u\partial v}dudv 
\ee
such that
$$\frac{\partial^2R_{H,K}(t,s)}{\partial t\partial s}= \alpha_{H,K}(t^{2H}+s^{2H})^{K-2}(ts)^{2H-1}-\beta_{H,K}|t-s|^{2HK-2}$$
where
$$\alpha_{H,K}=2^{-K+2}H^2K(K-1)$$
and
$$\beta_{H,K}=2^{-K+1}HK(2HK-1).$$
\vsp
If $\varphi, \psi \in |\clh_{W^{H,K}}|$, then their inner product in $\clh_{W^{H,K}}$ is given by
$$<\varphi,\psi>=\int_0^T\int_0^T \varphi(u) \psi(v)\frac{\partial^2R_{H,K}(u,v)}{\partial u\partial v}dudv.$$
Furthermore, for $\varphi, \psi \in |\clh_{W^{H,K}}|$, it can be checked that
$$E(\int_0^T\varphi(u)dW_u^{H,K})=0$$
and
$$E(\int_0^T\varphi(u)dW_u^{H,K}\int_0^T\psi(v)dW_v^{H,K})= <\varphi,\psi>_{\clh_{W^{H,K}}}.$$
The canonical Hilbert space $\clh_{W^{H,K}}$ associated with the bifBm $W^{H,K}$ satisfies the property:
$$L^2([0,T]) \subset L^{1/HK}([0,T])\subset |\clh_{W^{H,K}}| \subset \clh_{W^{H,K}}$$
when $H\in (0,1), K\in (0,1]$ such that $HK\in (\frac{1}{2},1).$
\vsp
\newsection{Preliminaries}
Let $W^{H,K}=\{W_t^{H,K}, t \geq 0\}$ be a bifractional Brownian motion with known parameters $H,K$ such that $H\in (0,1), K\in (0,1]$ and $HK\in (\frac{1}{2},1).$
Consider the problem of estimating the unknown function $\theta(t), 0 \leq t \leq T$ (linear multiplier) from the observations $\{X_t,0\leq t \leq T\}$ of process satisfying the stochastic differential equation
\be
dX_t=\theta(t) X_t dt + \epsilon \;dW_t^{H,K}, X_0=x_0, 0 \leq t \leq T
\ee
and study the properties of the estimator as $\epsilon \rightarrow 0.$ Consider the differential equation in the limiting system of (3.1), that is, for $\epsilon=0,$ given by
\be
dx_t=\theta(t) x_t dt, x_0, 0 \leq t \leq T.
\ee
Observe that
$$x_t=x_0 \exp \{ \int^t_0 \theta(s) ds).$$
\vsp
We assume that the following condition holds:
\vsp
\NI{$(A_1)$} The trend coefficient $\theta(t),$ over the interval $[0,T],$ is bounded by a constant $L$.
\vsp
\noindent{\bf Lemma 3.1.} Let the condition $(A_1)$ hold and  $\{X_t, 0\leq t \leq T\}$ and $\{x_t, 0\leq t \leq T\}$ be the solutions of the equations (3.1) and (3.2) respectively. Then, with probability one,
\be
|X_t-x_t| < e^{L t} \epsilon |W_t^{H,K}|
\ee
and
\be
\sup_{0 \leq t \leq T} E(X_t-x_t)^2  \leq e^{2L T} \epsilon^2 T^{2HK}.
\ee
\vsp
\noindent{\bf Proof of (a):} Let $u_t=|X_t-x_t|.$ Then by $(A_1)$; we have,
 \bea
u_t & \leq &\int^t_0 |\theta(v) (X_v-x_v)| dv + \epsilon |W_t^{H,K}|\\\nonumber
 & \leq & L \int^t_0 u_v dv + \epsilon |W_t^{H,K}|.\\\nonumber
\eea
Applying the Gronwall's lemma (cf. Lemma 1.12, Kutoyants (1994), p. 26), it follows that
\be
u_t \leq \epsilon |W_t^{H,K}| e^{L t}.
\ee
\vsp
\noindent{\bf Proof of (b):} From the equation (3.3), we have
\bea
E(X_t-x_t)^2 & \leq & e^{2 L t} \epsilon^2 E(|W_t^{H,K}|)^2\\\nonumber
&= & e^{2 L t} \epsilon^2 t^{2HK} . \\\nonumber
\eea
Hence
\be
\sup_{0 \leq t \leq T} E (X_t-x_t)^2 \leq e^{2 L T} \epsilon^2 T^{2HK}.
\ee
\vsp
\newsection{Main Results}
Let $\Theta_0(L)$ denote the class of all functions $\theta(.)$ with the same bound $L$. Let
$\Theta_k(L) $ denote the class of all functions $\theta(.)$ which are uniformly bounded by the same constant $L$
and which are $k$-times differentiable with respect to $t$ satisfying the condition
$$|\theta^{(k)}(x)-\theta^{(k)}(y)|\leq L_1|x-y|, x,y \in R$$
for some constant $L_1 >0.$ Here $g^{(k)}(x)$ denotes the $k$-th derivative of $g(.)$ at $x$ for $k \geq 0.$ If $k=0,$ we interpret the function $g^{(0)}(x)$ as $g(x).$
\vsp
Let $G(u)$ \ be a  bounded function  with compact support $[A,B]$ with $A<0<B$ satisfying the condition\\
\NI{$(A_2)$} $\int^B_A G(u) du =1.$
\vsp
It is obvious that the following conditions are satisfied by the function $G(.):$
\begin{description}
\item{(i)} $ \int^\infty_{-\infty} |G(u)|^2 du < \infty;$ \\
\item{(ii)}$\int^\infty_{-\infty} |u^{k+1} G(u)|^2 du <\infty.$\\
\end{description}
We define a kernel type estimator $\hat \theta_t$ of the function $\theta(t)$ by the relation
\be
\widehat{\theta}_t X_t= \frac{1}{\varphi_\epsilon}\int^T_0 G \left(\frac{\tau-t}{\varphi_\epsilon} \right) d X_\tau
\ee
where the normalizing function  $ \varphi_\epsilon \rightarrow 0 $ as \ $ \epsilon \rightarrow 0. $ Let $E_\theta(.)$ denote the expectation when the function $\theta(.)$ is the linear multiplier.
\vsp
\NI{\bf Theorem 4.1:}  {\it Suppose that the linear multiplier $\theta(.) \in \Theta_0(L)$ and  the function \ $ \varphi_\epsilon \rightarrow 0$  and $\epsilon^2\varphi_\epsilon^{2HK-2}\raro 0$ as $\epsilon \raro 0.$ Suppose the conditions $(A_1)-(A_2)$ hold.Then, for any $ 0 < a \leq b < T,$ the estimator $\hat \theta_t$ is uniformly consistent, that is,}
\be
\lim_{\epsilon \rightarrow 0} \sup_{\theta(.) \in \Theta_0(L)} \sup_{a\leq t \leq b } E_\theta ( |\hat \theta_t X_t-\theta(t) x_t|^2)= 0.
\ee
\vsp
In addition to the conditions $(A_1)$ and $(A_2),$ suppose the following condition holds:
\vsp
\NI{$(A_3)$}$ \int^\infty_{-\infty} u^j G(u)   du = 0 \;\;\mbox{for}\;\; j=1,2,...k.$
\vsp
\NI{\bf Theorem 4.2:} {\it Suppose that the function $ \theta(.) \in \Theta_{k+1}(L)$ and  the conditions $(A_1)-(A_3)$ hold. Further suppose that $
\varphi_\epsilon = \epsilon^{\frac{1}{k-HK+2}}.$ Then,
\be
\limsup_{\epsilon \rightarrow 0} \sup_{\theta(.) \in \Theta_{k+1}(L)}\sup_{a \leq t \leq b} E_\theta (| \hat\theta_t X_t - \theta(t)x_t|^2)
\epsilon^{-\min(2, \frac{2(k+1)}{k+2-HK})} \ < \infty.
\ee
\vsp
\NI{\bf Theorem 4.3:} {\it Suppose that the function $\theta(.) \in \Theta_{k+1}(L)$ for some $k>1$ and  the conditions $(A_1)-(A_3)$ hold. Further suppose that 
$\varphi_\epsilon= \epsilon^{\frac{1}{k-HK+2}}.$ Let $J(t)= \theta(t)x_t.$ Then, as $\epsilon \raro 0,$  the asymptotic distribution of
$$ \epsilon^{\frac{-(k+1)}{k-HK+2}} (\hat\theta_t X_t - J(t)  - \frac{J^{(k+1)}(t)}{(k+1) !} \int^\infty_{-\infty} G(u) u^{k+1}\ du)$$
is  Gaussian mean zero and variance
$$ \sigma^2_{H,K}= \int^{\infty}_{-\infty} \int^{\infty}_{-\infty}G(u)G(v) \frac{\partial^2 R_{H,K}(u,v)}{\partial u\partial v}\ dudv$$
where
$$\frac{\partial^2R_{H,K}(t,s)}{\partial t\partial s}= \alpha_{H,K}(t^{2H}+s^{2H})^{K-2}(ts)^{2H-1}-\beta_{H,K}|t-s|^{2HK-2},$$

$$\alpha_{H,K}=2^{-K+2}H^2K(K-1),$$
and}
$$\beta_{H,K}=2^{-K+1}HK(2HK-1).$$}
\vsp
\newsection{Proofs of Theorems}
\NI{\bf Proof of Theorem 4.1 :} From the inequality
$$(a+b+c)^2\leq 3(a^2+b^2+c^2), a,b,c\in R,$$
it follows that
\bea
\;\;\;\\\nonumber
E_\theta[|\hat \theta (t) x_t -\theta(t) x_t|^2]  &=& E_\theta [ |\frac{1}{\varphi_\epsilon}  \int^T_0 G \left(\frac{\tau-t}{\varphi_\epsilon} \right) \left(\theta(\tau) X_\tau -\theta(\tau) x_\tau \right)  d \tau \\ \nonumber
& &+ \frac{1}{\varphi_\epsilon} \int^T_0 G \left(\frac{\tau-t}{\varphi_\epsilon}\right) \theta(\tau) x_\tau d \tau- \theta(t) x_t
 + \frac{\epsilon}{\varphi_\epsilon} \int^T_0 G \left(\frac{\tau-t}{\varphi_\epsilon} \right)
 dW_\tau^{H,K}|^2]\\ \nonumber
 & \leq  & 3 E_\theta[ |\frac{1}{\varphi_\epsilon}  \int^T_0 G \left(\frac{\tau-t}{\varphi_\epsilon} \right) (\theta(\tau) X_\tau -\theta(\tau) x_\tau) d\tau|^2]\\ \nonumber
 & & + 3 E_\theta [|\frac{1}{\varphi_\epsilon} \int^T_0 G \left(\frac{\tau-t}{\varphi_\epsilon} \right)\theta(\tau) x_\tau d\tau -\theta(t)x_t |^2 ]\\ \nonumber
 & & +  3 \frac{\epsilon^2}{\varphi_\epsilon^2} E_\theta [ |\int^T_0 G \left( \frac{\tau-t}{\varphi_\epsilon}\right) dW_\tau^{H,K}|^2]\\ \nonumber
 &= & I_1+I_2+I_3 \;\;\mbox{(say).}\;\;\\ \nonumber
\eea
By the boundedness condition on the function $\theta(.),$ the inequality (3.3) in Lemma 3.1 and the condition $(A_2)$, and applying the H\"older inequality, it follows that
\bea
\;\;\;\\\nonumber
I_1 &= &3 E_\theta \left| \frac{1}{\varphi_\epsilon} \int^T_0 G
\left(\frac{\tau-t}{\varphi_\epsilon} \right) (\theta(\tau) X_\tau -\theta(\tau) x_\tau)
d\tau \right|^2 \\\nonumber
&= & 3E_\theta  \left| \int^\infty_{-\infty} G(u) \left(\theta(t+\varphi_\epsilon u) X_{t+\varphi_\epsilon u}  - \theta(t+\varphi_\epsilon u) x_{t+\varphi_\epsilon u}\right) du\right|^2\\\nonumber
& \leq & 3 (B-A) \int^\infty_{-\infty} |G(u)|^2 L^2 E \left|X_{t+\varphi_\epsilon u}-x_{t+\varphi_\epsilon u} \right|^2 \ du
\;\;\mbox{(by using the condition $(A_1)$)}\\\nonumber
& \leq & 3(B-A)\int^\infty_{-\infty} |G(u)|^2 \;\;L^2 \sup_{0 \leq t +
\varphi_\epsilon u \leq T}E_\theta \left|X_{t+\varphi_\epsilon u}
-x_{t+\varphi_\epsilon u}\right|^2 \ du \\\nonumber
& \leq & 3 (B-A)L^2  e^{2LT} \epsilon^2 T^{2HK} \int_{-\ity}^\ity|G(u)|^2du\;\;\mbox{(by using (3.4))}\\\nonumber
\eea
which tends to zero as $\epsilon \raro 0.$  For the term $I_2$, by the boundedness condition on the function $\theta(.),$ the condition $(A_2)$ and the H\"older inequality, it follows that
\bea
\;\;\;\\\nonumber
I_2 &= & 3E_\theta \left| \frac{1}{\varphi_\epsilon} \int^T_0 G\left(
\frac{\tau-t}{\varphi_\epsilon}\right) \theta(\tau) x_\tau d \tau - \theta(t) x_t\right|^2 \\ \nonumber
& = & 3  \left| \int^\infty_{-\infty} G(u)
\left(\theta(t+\varphi_\epsilon u) x_{t+\varphi_\epsilon u}-\theta(t) x_t \right)  \ du \right|^2
\\\nonumber
& \leq  & 3 (B-A)L^2 \varphi_\epsilon^2  \int_{-\ity}^\ity|uG(u)|^2 du\;\;\mbox{(by $(A_2)$)}.\\\nonumber
\eea
The last term  tends to zero as  $\varphi_\epsilon \rightarrow 0.$ We will now get an upper bound on the term $I_3.$ Note that
\bea
\;\;\;\\ \nonumber
I_3 &= & 3\frac{ \epsilon^2}{\varphi_\epsilon^2} E_\theta \left|\int^T_0 G \left(\frac{\tau-t}{\varphi_\epsilon}\right ) dW_\tau^{H,K}\right|^2 \\ \nonumber
&=& 3 \frac{ \epsilon^2}{\varphi_\epsilon^2} \int_0^T\int_0^T G \left(\frac{\tau-t}{\varphi_\epsilon}\right )G \left(\frac{\tau^\prime-t}{\varphi_\epsilon}\right ) \frac{\partial^2 R_{H,K}(s,t)}{\partial s\partial t}|_{s=\tau,s^\prime=\tau^\prime} d\tau d\tau^\prime\\\nonumber
& \leq & C_1\frac{\epsilon^2}{\varphi^2_\epsilon}\varphi^2_\epsilon \int_R\int_R G(u) G(v)[\alpha_{H,K}((t-\varphi_\epsilon u)^{2H} + (t-\varphi_\epsilon v)^{2H})^{K-2})((t-\varphi_\epsilon u)(t-\varphi_\epsilon v))^{2H-1}\\\nonumber
&&-\beta_{H,K}|(t-\varphi_\epsilon u)-(t-\varphi_\epsilon v)|^{2HK-2}]dudv\\\nonumber
&\leq & C_2\epsilon^2[t^{2H(K-2)+2(2H-1)}+(\varphi_\epsilon)^{2HK-2}]\\\nonumber
&\leq & C_3\epsilon^2+ C_3\epsilon^2(\varphi_\epsilon)^{2HK-2}.
\eea
for some positive constant $C_3.$ Theorem 4.1 is now proved by using the equations (5.1) to (5.4).
\vsp
\NI {\bf Proof of Theorem 4.2 :} Let $J(t)=\theta(t) x_t.$ By the Taylor's formula, for any $x \in R,$
$$ J(y) = J(x) +\sum^k_{j=1} J^{(j)} (x) \frac{(y-x)^j}{j !} +[ J^{(k)} (z)-J^{(k)} (x)] \frac{(y-x)^k}{k!} $$
for some $z$ such that $|z-x|\leq |y-x|.$ Using this expansion, the equation (3.2) and the condition $(A_3)$  in the expression for $I_2$ defined in the proof of  Theorem 4.1, it follows that
\ben
\;\;\\\nonumber
I_2 & = & 3 \left[
\int^\infty_{-\infty} G(u) \left(J(t+\varphi_\epsilon u) - J(t) \right)  \ du \right]^2\\ \nonumber
&= & 3[ \sum^k_{j=1} J^{(j)} (t) (\int^\infty_{-\infty}G(u) u^j du )\varphi^j_\epsilon (j!)^{-1}\\\nonumber
& & \;\;\;\;+(\int^\infty_{-\infty}G(u) u^k (J^{(k)}(z_u) -J^{(k)} (x_t))du \;\varphi^k_\epsilon (k !)^{-1}]^2\\ \nonumber
\een
for some $z_u$ such that $|x_t-z_u|\leq |x_{t+\varphi_\epsilon u}-x_t| \leq C|\varphi_\epsilon u|.$ Hence
\bea
I_2 & \leq & 3  L^2 \left[  \int^\infty_{-\infty} |G(u)u^{k+1}|\varphi^{k+1}_\epsilon (k!) ^{-1}  du  \right]^2
\\ \nonumber
& \leq & 3 L^2 (B-A)(k!)^{-2} \varphi^{2(k+1)}_\epsilon\int^\infty_{-\infty} G^2(u) u^{2 (k+1)}\ du \\\nonumber 
&\leq & C_2 \varphi_\epsilon^{2(k+1)}\\ \nonumber
\eea
for some positive constant $C_2$. Combining the equations (5.2)- (5.5), we get that there exists a positive constant $C_3$
such that 
$$ \sup_{a \leq t \leq b}E_\theta|\hat\theta_t X_t-\theta(t) x_t|^2 \leq C_3 (\epsilon^2 +  \varphi^{2(k+1)}_\epsilon +\epsilon^2 \varphi_\epsilon^{2HK-2}). $$ 
Choosing $ \varphi_\epsilon = \epsilon^{\frac{1}{k+2-HK}},$  we get that 
$$ \limsup_{\epsilon\rightarrow 0} \sup_{\theta(.) \in \Theta_{k+1} (L) } \sup_{a \leq t
\leq b} E_\theta|\theta(t)X_t - \theta(t)x_t|^2\epsilon^ {-\min(2, \frac{2(k+1)}{k+2-HK})} < \infty. $$ This completes the proof of
Theorem 4.2. 
\vsp 
\NI{\bf Proof of Theorem 4.3:} Let $\alpha= \frac{k+1}{k-HK+2}.$ Note that $0<\alpha<1$ since $0<HK<1.$ From (3.1), we obtain that
\bea
\lefteqn{\epsilon^{-\alpha}( \hat\theta(t)X_t -\theta(t)x_t)}\\\nonumber
 &= &\epsilon^{-\alpha}[\frac{1}{\varphi_\epsilon} \int^T_0 G \left(\frac{\tau-t}{\varphi_\epsilon} \right)
 \left( \theta(\tau)X_\tau-\theta(\tau) x_\tau\right) \  d \tau \\ \nonumber
 & & + \frac{1}{\varphi_\epsilon} \int^T_0 G \left( \frac{\tau-t}{\varphi_\epsilon}\right) \theta(\tau) x_\tau d\tau -\theta(t) x_t+ \frac{\epsilon}{\varphi_\epsilon} \int^T_0 G \left( \frac{\tau-t}{\varphi_\epsilon}\right) dW_\tau^{H,K}]\\ \nonumber
 &= & \epsilon^{-\alpha} [ \int^\infty_{-\infty} G(u) (\theta(t+\varphi_\epsilon u)X_{t+\varphi_\epsilon u} - \theta(t+\varphi_\epsilon u) x_{t+\varphi_\epsilon u}) \ du  \\ \nonumber
 & & +\int^\infty_{-\infty} G(u) ( \theta(t+\varphi_\epsilon u) x_{t+\varphi_\epsilon u}- \theta(t) x_t) \ du \\ \nonumber
 & &+ \frac{\epsilon}{\varphi_{\epsilon}}\int^T_0 G \left(\frac{\tau-t}{\varphi_\epsilon} \right) dW_\tau^{H,K}].\\ \nonumber
&=& R_1+R_2+R_3 \;\;\;\mbox{(say).}\\\nonumber
\eea
By the boundedness  condition on the function $\theta(.)$ and part (a) of Lemma 3.1, it follows that
\bea
R_1 & \leq & \epsilon^{-\alpha}|\int_{-\ity}^\ity G(u)(\theta(t+\varphi_\epsilon u) X_{t+\varphi_\epsilon u} - \theta(t+\varphi_\epsilon u)x_{t+\varphi_\epsilon u}) du|\\\nonumber
&\leq & \epsilon^{-\alpha} \epsilon L \int_{-\ity}^\ity |G(u)| X_{t+\varphi_\epsilon u} - x_{t+\varphi_\epsilon u}| du\\\nonumber
& \leq & Le^{LT} \epsilon^{1-\alpha} \int_{-\ity}^\ity |G(u)|\sup_{0\leq t+\varphi_\epsilon u \leq T}|W_{t+\varphi_\epsilon u}^{H,K}|du.\\\nonumber
\eea
Applying the Markov's inequality, it follows that, for any $\eta >0,$
\bea
P(|R_1|>\eta) &\leq &  \epsilon^{1-\alpha} \eta^{-1} Le^{LT} \int_{-\ity}^\ity |G(u)|E_\theta(\sup_{0\leq t+\varphi_\epsilon u \leq T}|W^{H,K}_{t+\varphi_\epsilon u}|)du\\\nonumber
&\leq & \epsilon^{1-\alpha}\eta^{-1} Le^{LT} \int_{-\ity}^\ity |G(u)||E_\theta[(\sup_{0\leq t+\varphi_\epsilon u \leq T}(W^{H,K}_{t+\varphi_\epsilon u})^2]|^{1/2} du\\\nonumber
&\leq & \epsilon^{1-\alpha}  \eta^{-1} Le^{LT} C T^{HK}\int_{-\ity}^\ity|G(u)|du\\\nonumber
\eea
from the maximal inequality for a bifBm proved in Theorem 2.3 in Prakasa Rao (2024) for some constant $C>0,$ and the last term tends to zero as $\epsilon \raro 0.$ Let $J_t=\theta(t)x_t.$ By the Taylor's formula, for any $t \in [0,T],$
$$ J_t = J_{t_0} + \sum^{k+1}_{j=1} J_{t_0}^{(j)}  \frac{(t-t_0)^j}{j !} + [ J_{t_0+\gamma(t-t_0)}^{(k+1)}-J_{t_0}^{(k+1)}] \frac{(t-t_0)^{k+1}}{(k+1)!} $$
where $0<\gamma<1$ and $t_0 \in (0,T).$ Applying the Condition $(A_3)$ and the Taylor's expansion, it follows that
\bea
R_2 &=& \epsilon ^{-\alpha}[\sum_{j=1}^{k+1}J_t^{(j)}(\int_{-\ity}^\ity G(u) u^j \;du)\varphi_\epsilon^j(j!)^{-1}\\\nonumber
&& \;\;\;\; +\frac{\varphi_\epsilon^{k+1}}{(k+1)!}\int_{-\ity}^\ity G(u) u^{k+1}(J_{t+\gamma \varphi_\epsilon u}^{(k+1)}-J_t^{(k+1)})\;du]\\\nonumber
&=& \epsilon ^{-\alpha} \frac{J_t^{(k+1)}}{(k+1)!}\int_{-\ity}^\ity G(u) u^{k+1}\;du\\\nonumber
&& \;\;\;\; + \varphi_\epsilon^{k+1} \epsilon^{-\alpha}\frac{1}{(k+1)!}\int_{-\ity}^\ity G(u)u^{k+1}(J_{t+\gamma  \varphi_\epsilon u}^{(k+1)}-J_t^{(k+1)})\;du.\\\nonumber.
\eea
Observing that $\theta(t) \in \Theta_{k+1}(L),$ we obtain that
\bea
\lefteqn{\frac{1}{(k+1)!}\int_{-\ity}^\ity G(u)u^{k+1}(J_{t+\gamma\varphi_\epsilon u}^{(k+1)}-J_t^{(k+1)})du}\\\nonumber
&\leq & \frac{1}{(k+1)!}\int_{-\ity}^\ity |G(u)u^{k+1}(J_{t+\gamma\varphi_\epsilon u}^{(k+1)}-J_t^{(k+1)})|du\\\nonumber
&\leq & \frac{L\varphi_\epsilon}{(k+1)!}\int_{-\ity}^\ity |G(u)u^{k+2}|du.\\\nonumber
\eea
Combining the equations given above, it follows that
\bea
\lefteqn{\epsilon^{-\alpha} (\hat\theta_tX_t-J(t)- \frac{J_t^{(k+1)}}{(k+1)!} \int^\infty_{-\infty} G(u) u^{k+1}\ du)}\\\nonumber
&=& O_p(\epsilon^{1-\alpha})+O_p(\epsilon^{-\alpha}\varphi_\epsilon^{k+2})+\epsilon^{1-\alpha}\varphi_\epsilon^{-1}\int_0^TG(\frac{\tau-t}{\varphi_\epsilon})dW_\tau^{H,K}.
\eea
Let
\be
\eta_\epsilon(t)= \epsilon^{\frac{-(k+1)}{k-HK+2}}\epsilon\varphi_\epsilon^{-1}\int_0^TG(\frac{\tau-t}{\varphi_\epsilon})dW_\tau^{H,K}.
\ee
Note that $E[\eta_\epsilon(t)]=0,$ and 
\ben
E([\eta_\epsilon(t)]^2)&=& (\epsilon^{\frac{1-HK}{k-HK+2}}\varphi_\epsilon^{-1})^2E([\int_0^TG(\frac{\tau-t}{\varphi_\epsilon})dW_\tau^{H,K}]^2)\\
&=& (\epsilon^{\frac{1-HK}{k-HK+2}}\varphi_\epsilon^{-1})^2[\varphi_\epsilon^{2HK}\int_R\int_R G(u)G(v)\frac{\partial^2 R_{H,K}(u,v)}{\partial u \partial v}dudv].
\een
Choosing $\varphi_\epsilon=\epsilon^{\frac{1}{k-HK+2}},$ we get that
$$E([\eta_\epsilon(t)]^2)=\int_R\int_R G(u)G(v)\frac{\partial^2R_{H,K}(u,v)}{\partial u \partial v}dudv.$$

From the choice of $\varphi_\epsilon$ and $\alpha,$ it  follows that 
$$ \epsilon^{1-\alpha}\varphi_\epsilon^{-1}= \varphi_\epsilon^{HK}$$
and, 
\bea
\;\;\;\;\\\nonumber
\lefteqn{Var[ \varphi_\epsilon^{-HK} \int^T_0 G \left(\frac{\tau-t}{\varphi_\epsilon} \right) dW^{H,K}_\tau]}\\\nonumber
&=& \varphi_\epsilon^{-2HK}\int_0^T\int_0^TG\left(\frac{\tau-t}{\varphi_\epsilon} \right) G\left(\frac{\tau^\prime-t}{\varphi_\epsilon} \right) \frac{\partial^2 R_{H,K}}{\partial s \partial s^\prime}|_{s=\tau,s^\prime=\tau^\prime}d\tau d\tau^\prime\\\nonumber
\eea
and the last term tends to 
$$\int_R\int_R G(u) G(v) \frac{\partial^2R_{H,K}(u,v)}{\partial u\partial v} dudv= \sigma^2_{H,K}$$
as $\epsilon \raro 0.$
Applying the Slutsky's theorem and the equations derived above, it can be checked that the random variable 
$$\epsilon^{-\alpha} (\hat\theta_t X_t- J_t - \frac{J_t^{(k+1)} }{(k+1) !} \int^\infty_{-\infty} G(u) u^{k+1}\ du)$$
has a limiting distribution as $\epsilon \raro 0$ as that of the family of random variables 
$$\varphi_\epsilon^{-HK}\int_{-\infty}^{\infty} G\left(\frac{\tau-t}{\varphi_\epsilon} \right) dW^{H,K}_\tau$$
as $\epsilon \raro 0$ which is Gaussian with mean zero and variance $\sigma^2_{H,K}.$ This completes the proof of Theorem 4.3.
\vsp
\newsection{Alternate Estimator for the Multiplier $\theta(.)$}
Let $\Theta_\rho(L_\gamma)$ be a class of functions $\theta(t)$ uniformly bounded  by a constant $L$ and $k$-times continuously differentiable for some integer 
$k \geq 1 $ with the $k$-th derivative satisfying the H\"older condition of the order $\gamma \in (0,1):$
$$|\theta^{(k)}(t)-\theta^{(k)}(s)|\leq L_\gamma |t-s|^\gamma, \rho=k+\gamma$$
and suppose that $\rho>HK.$ Suppose the process $\{X_t,0\leq t \leq T\}$ satisfies the stochastic differential equation given by the equation (3.1) where the linear multiplier is an unknown function in the class $\Theta_\rho(L_\gamma)$ and further suppose that $x_0 > 0$ and is {\it known}. 
From the Lemma 3.1, it follows that
$$|X_t-x_t| \leq \epsilon e^{Lt}\sup_{0\leq s \leq T} |W^{H,K}_s|.$$
Let
$$A_t= \{\omega: \inf_{0\leq s \leq t}X_s(\omega)\geq \frac{1}{2}x_0e^{-Lt}\}$$
and let $A=A_T.$ Following the technique suggested in Kutoyants (1994), p. 156, we define another process $Y$ with the differential
$$dY_t=\theta(t) I(A_t) dt + \epsilon 2x_0^{-1} e^{LT} I(A_t)\;dW^{H,K}_t, 0\leq t \leq T.$$
We will now construct an alternate estimator of the linear multiplier $\theta(.)$ based on the process $Y$ over the interval $[0,T].$
Define the estimator
$$\tilde \theta(t)= I(A) \frac{1}{\varphi_\epsilon}\int_0^T G(\frac{t-s}{\varphi_\epsilon})dY_s$$
where the kernel function $G(.)$ satisfies the conditions $(A_1)-(A_3)$. Observe that
\ben
E|\tilde \theta(t)-\theta(t)|^2 &= & E_\theta|I(A) \frac{1}{\varphi_\epsilon}\int_0^T G(\frac{t-s}{\varphi_\epsilon})(\theta(s)-\theta(t))ds\\\nonumber
&& \;\;\;\; + I(A^c)\theta(t)+I(A)\frac{\epsilon}{\varphi_\epsilon}\int_0^T G(\frac{t-s}{\varphi_\epsilon})2x_0^{-1}e^{LT}dW^{H,K}_s|^2\\\nonumber
&\leq & 3 E_\theta|I(A)\int_R G(u)[\theta(t+u\varphi_\epsilon)-\theta(t)]du|^2+ 3 |\theta(t)|^2 [P(A^c)]^2\\\nonumber
&& \;\;\;\; + 3 \frac{\epsilon^2}{\varphi_\epsilon^2}|E[I(A)\int_0^T G(\frac{t-s}{\varphi_\epsilon})2x_0^{-1}e^{LT}dW^{H,K}_s]|^2\\\nonumber
&=& D_1+D_2+D_3. \;\;\mbox{(say)}.\\\nonumber
\een
Applying the Taylor's theorem and using the fact that the function $\theta(t)\in \Theta_\rho(L_\gamma)$, it follows that
\ben
D_1 \leq C_1\frac{1}{(k+1)!}\varphi_\epsilon^{2\rho} \int_R|G^2(u)u^{2\rho}|du.
\een
Note that, by Lemma 3.1,
\ben
P(A^c) & = & P(\inf_{0\leq t \leq T}X_t < \frac{1}{2}x_0e^{-LT})\\\nonumber
&\leq & P(\inf_{0 \leq t \leq T}|X_t-x_t| + \inf_{0\leq t \leq T}x_t < \frac{1}{2}x_0e^{-LT})\\\nonumber
&\leq & P(\inf_{0 \leq t \leq T}|X_t-x_t| < -\frac{1}{2}x_0e^{-LT})\\\nonumber
&\leq & P(\sup_{0 \leq t \leq T}|X_t-x_t| > \frac{1}{2}x_0e^{-LT})\\\nonumber
&\leq & P(\epsilon e^{LT}\sup_{0 \leq t \leq T}|W_t^{H,K}|>\frac{1}{2}x_0e^{-LT})\\\nonumber
&= & P(\sup_{0 \leq t \leq T}|W_t^{H,K}|>\frac{x_0}{2\epsilon}e^{-2LT})\\\nonumber
&\leq & (\frac{x_0}{2\epsilon}e^{-2LT})^{-2}E[\sup_{0 \leq t \leq T}|W^{H,K}_t|^2]\\\nonumber
&\leq & (\frac{x_0}{2\epsilon}e^{-2LT})^{-2}C_2T^{2HK}
\een
by Theorem 2.3 in Prakasa Rao (2024) for some positive constant $C_2.$ The upper bound obtained above and the fact that $|\theta(s)|\leq L, 0\leq s \leq T$ leads  an upper bound for the term $D_2.$ We have used the inequality
$$x_t= x_0 \exp(\int_0^t\theta (s)ds)\geq x_0 e^{-Lt}$$
in the computations given above. Applying Theorem 2.1, it follows that
\ben
\lefteqn{E[|I(A)\int_0^T G(\frac{t-s}{\varphi_\epsilon})2x_0^{-1}e^{LT}dW_s^{H,K}|^2]}\\\nonumber
&\leq & CE[|\int_0^T G(\frac{t-s}{\varphi_\epsilon})dW^{H,K}_s|^2]\\\nonumber
&=& C\;Var[\int_0^T G(\frac{t-s}{\varphi_\epsilon})dW^{H,K}_s]\\\nonumber
&=& C\varphi_\epsilon^{2HK}\int_R\int_RG(u)G(v) \frac{\partial^2R_{H,K}(u,v)}{\partial u \partial v}dudv\\\nonumber
\een
for some positive constant $C$ which leads to an upper bound on the term $D_3.$ Combining the above estimates, it follows that
\ben
E|\tilde \theta(t)-\theta(t)|^2\leq C_1\varphi_\epsilon^{2\rho} + C_2 \epsilon^4+ C_3\epsilon^2 \varphi_\epsilon^{2HK}
\een
for some positive constants $C_i, i=1,2,3.$ Choosing $\varphi_\epsilon=\epsilon^{\frac{1}{\rho-HK}},$ we obtain that
\ben
E|\tilde \theta(t)-\theta(t)|^2\leq C_4 \epsilon^{\frac{2\rho}{\rho-HK}}+ C_5 \epsilon^{4}
\een
for some positive constants $C_4$ and $C_5.$  Hence we obtain the following result implying the uniform consistency of the estimator $\tilde \theta(t)$ as an estimator of $\theta(t)$ as $\epsilon \raro 0.$ 
\vsp
\NI{\bf Theorem 6.1:} {\it Let $\theta \in \Theta_\rho(L)$ where $\rho >HK.$ Let $\varphi_\epsilon= \epsilon^{1/(\rho-HK)}.$ Suppose the conditions $(A_1)-(A_3)$ hold. Then, for any interval $[a,b] \subset [0,T],$ }
\ben
\limsup_{\epsilon \raro 0}\sup_{\theta(.)\in \Theta_\rho(L)}\sup_{a\leq t \leq b}E|\tilde \theta(t)-\theta(t)|^2 \epsilon^{-\min(4, \frac{2\rho}{\rho-HK})}<\ity.
\een
\vsp

\NI{\bf Acknowledgment:} This work was supported by the Indian National Science Academy (INSA) under the scheme  ``INSA Honorary Scientist" at the CR RAO Advanced Institute of Mathematics, Statistics and Computer Science, Hyderabad 500046, India.
\vsp
\NI{\bf References :}
\vsp
\begin{description}

\item Bardina, X. and Es-Sebaiy, K. 2011. An extension of bifractional Brownian motion, {\it communications on Stochastic Analysis}, 5: 333-340.

\item Houdr\'e, C. and Villa, J. 2003. An example of infinite dimensional quasi-helix, {\it Contemp. Math., Amer. Math. Soc.}, 336: 195-201.

\item Kahane, J.P. 1981. H\'elices et quasi-h\'elices, {\it Adv.Math.} B 7, 417-433.

\item Kahane, J.P. 1985. {\it Some Random Series Functions}, Cambridge University Press, Cambridge.

\item Keddi, A., Madani, F. and Bouchentouf, A.A. 2020. Nonparametric estimation of trend function for stochastic differential equations driven by a bifractional Brownian motion, {\it Acta Univ. Sapientiae, Mathematica}, 12: 128-145.

\item Kruk, I., Russo, F. and Tudor, C.A. 2007. Wiener integrals, Malliavin calculus and covariance measure structure, {\it J. Funct. Anal. }, 249: 92-142.
    
\item Kutoyants, Y.A. 1994. {\it Identification of Dynamical Systems with small Noise}, Dordrecht, Kluwer. 

\item Lipschits, M. and Volkova, K. 2015. Bifractional Brownian motion: Existence and border cases, arXiv:1502.02217.    

\item Mishra, M.N. and Prakasa Rao, B.L.S. 2011. Nonparametric estimation of linear multiplier for fractional diffusion processes, {\it Stochastic Anal. Appl.}, 29: 706-712.

\item Mishura, Y. 2008. {\it Stochastic Calculus for Fractional Brownian Motion and Related Processes}, Berlin, Springer.

\item Prakasa Rao, B.L.S. 2010. {\it Statistical Inference for Fractional Diffusion Processes}, London, Wiley.

\item Prakasa Rao, B.L.S. 2024. Maximal inequalities for bifractional Brownian motion, Preprint, CR Rao AIMSCS, Hyderabad, India.

\item Russo, F. and Tudor, C. 2006. On the bifractional Brownian motion, {\it Stoch. Process. Their Appl.}, 116: 830-856.

\item Tudor, C. 2013.  {\it Analysis of Variations for Self-similar Processes}, Springer, Switzerland.

\item Tudor, C. 2023. {\it Stochastic Partial Differential Equations with Additive Gaussian Noise}, World Scientific, Singapore.

\item Tudor, C. and Xiao, Y. 2007. Sample path properties of bifractional Brownian motion, {\it Bernoulli}, 13: 1023-1052.

\end{description}

\end{document}